\documentclass[11pt]{article}
\usepackage{amsthm}
\usepackage{latexsym, amssymb, amsmath}
\usepackage{graphicx}
\usepackage[dvipsnames]{xcolor}
\usepackage{fancyvrb}

\renewcommand{\theequation}{\arabic{section}.\arabic{equation}} 
\begin{document}
\numberwithin{equation}{section}

\newtheorem{THEOREM}{Theorem}
\newtheorem{PRO}{Proposition}
\newtheorem{XXXX}{\underline{Theorem}}
\newtheorem{CLAIM}{Claim}
\newtheorem{COR}{Corollary}
\newtheorem{LEMMA}{Lemma}
\newtheorem{REM}{Remark}
\newtheorem{EX}{Example}
\newenvironment{PROOF}{{\bf Proof}.}{{\ \vrule height7pt width4pt depth1pt} \par \vspace{2ex} }
\newcommand{\Bibitem}[1]{\bibitem{#1} \ifnum\thelabelflag=1 
  \marginpar{\vspace{0.6\baselineskip}\hspace{-1.08\textwidth}\fbox{\rm#1}}
  \fi}
\newcounter{labelflag} \setcounter{labelflag}{0}
\newcommand{\labelon}{\setcounter{labelflag}{1}}
\newcommand{\Label}[1]{\label{#1} \ifnum\thelabelflag=1 
  \ifmmode  \makebox[0in][l]{\qquad\fbox{\rm#1}}
  \else\marginpar{\vspace{0.7\baselineskip}\hspace{-1.15\textwidth}\fbox{\rm#1}}
  \fi \fi}

\newcommand{\LEFTLINE}{\ifhmode\newline\else\noindent\fi}
\newcommand{\RIGHTLINE}[1]{\LEFTLINE\rightline{#1}}
\newcommand{\CENTERLINE}[1]{\LEFTLINE\centerline{#1}}
\def\BOX #1 #2 {\framebox[#1in]{\parbox{#1in}{\vspace{#2in}}}}
\parskip=8pt plus 2pt
\def\AUTHOR#1{\author{#1} \maketitle}
\def\Title#1{\begin{center}  \Large\bf #1 \end{center}  \vskip 1ex }
\def\Author#1{\vspace*{-2ex}\begin{center} #1 \end{center}  
 \vskip 2ex \par}
\renewcommand{\theequation}{\arabic{section}.\arabic{equation}}
\def\bdk#1{\makebox[0pt][l]{#1}\hspace*{0.03ex}\makebox[0pt][l]{#1}\hspace*{0.03ex}\makebox[0pt][l]{#1}\hspace*{0.03ex}\makebox[0pt][l]{#1}\mbox{#1} }
\def\psbx#1 #2 {\mbox{\psfig{file=#1,height=#2}}}

 
\newcommand{\ia}{\,\,\nearrow}
\newcommand{\da}{\,\,\searrow}
\newcommand{\ic}{\nearrow}
\newcommand{\dc}{\searrow}
\newcommand{\me}{\mbox{e}}
\newcommand{\vp}[2]{\vphantom{\vrule height#1pt depth#2pt}}
\newcommand{\FG}[2]{{\includegraphics[width=#1mm]{#2.eps}}}
\renewcommand{\thefootnote}{\fnsymbol{footnote}}
\newcommand{\sr}{{\par\centerline{\vrule height0.1pt width2.5in depth0pt}\par}}
 
 
\renewcommand{\theequation}{\arabic{equation}}
\Title{A Numerical Procedure for Proving Specific Strict One-Variable Inequalities in Specific Finite Intervals}


\vspace{0.5cm}
\begin{center}
MAN KAM KWONG
\end{center}

\begin{center}
\emph{Department of Applied Mathematics\\ The Hong Kong Polytechnic University,\\ Hunghom, Hong Kong}\\
\tt{mankwong@polyu.edu.hk}
\end{center}

\par\vspace*{\baselineskip}\par

\newcommand{\Cr}{\color{red}}
\newcommand{\mb}{\mathbf}
\newcommand{\tx}[1]{\mbox{ #1 }}

\parskip=6pt

\begin{abstract}
A numerical procedure and its MAPLE implementation capable of rigorously,
albeit in a brute-force manner, proving
specific strict one-variable inequalities in specific finite intervals
is described. The procedure is useful, for instance, to affirm strict lower
bounds of specific functions.
\end{abstract}

\vspace{0.9cm}

\vspace{0.2cm}

\par\vspace*{\baselineskip}\par
{\par\vspace*{\baselineskip}\hrule{\vspace*{\baselineskip}\par}}

\section{Introduction }
Is there a good way to prove 
{\small 
\begin{equation}  \left( \frac{1-3x}{2} \right) \ln\left( \frac{1-3x}{2} \right)  + 2 \left( \frac{1-24x}{5} \right) \ln\left( \frac{1-24x}{5} \right)  > \,\, 3 \left( \frac{1-15x}{4} \right) \ln\left( \frac{1-15x}{4} \right)   \Label{ln}  \end{equation}
}
for $ x\in[0,0.04] $, or
{\small 
\begin{equation}  - \frac{\ln^3(2)}{2^s} + \frac{\ln^3(3)}{3^s}  -  \frac{\ln^3(4)}{4^s} + \frac{\ln^3(5)}{2\cdot5^s}  > 0   \Label{b4p}  \end{equation}
}
for $ s\in[0,1] $?

Each is a strict inequality involving only one variable (no parameters) in a 
specific finite interval.
Such inequalities are rarely of general interest. No one will present
any of them as a lemma, much less as a theorem. Yet, often in the proof of
some general result, the need for verifying an inequality like that pops up. For instance,
(\ref{b4p}) appears in \cite{ak1} and (\ref{ln}) in the proof of the non-negativity 
of a class of sine polynomials in $ [0,\pi ] $.

Numerical computation and graphing
provide evidence for the validity of such an
inequality, but obviously not as rigorous proofs.
Nevertheless, a simple principle makes it possible to use numerics
to rigorously prove the inequality. I dub it the DIF technique.

\begin{center}
\fbox{
\begin{minipage}{5in}
\par\vspace*{1mm}\par
\em Let $ f(x)=g_1(x)-g_2(x) $ be the difference of two increasing (decreasing) functions in an interval
$ [\alpha ,\beta ] $. To prove that $ f(x)>0 $ in $ [\alpha ,\beta ] $, it suffices to exhibit a
sequence of points $ \left\{ \tau _k\right\} _{k=1}^n $ such that
$$  \alpha  = \tau _1 < \tau _2 < \cdots < \tau _n = \beta   $$
and
\par\vspace*{-6mm}\par
\begin{equation}  g_2(\tau _{i+1}) < g_1(\tau _{i}), \qquad  \forall i = 1,2,...,n-1 .  \Label{g12}  \end{equation}
$$  \big(\,  g_2(\tau _{i}) < g_1(\tau _{i+1})  \,\big) \hspace*{33mm}  $$
\par\vspace*{0.0mm}\par
\end{minipage}
}
\end{center}

The obvious proof is omitted. The requirement of having to verify the inequality for
all $ x\in[\alpha ,\beta ] $ reduces to checking (\ref{g12}) only 
for a finite number of values. 
   
In applying the principle to (\ref{ln}), the increasing functions
$$  g_1(x) =  \left( \frac{1-3x}{2} \right) \ln\left( \frac{1-3x}{2} \right)  + 2 \left( \frac{1-24x}{5} \right) \ln\left( \frac{1-24x}{5} \right)    $$
and 
$$  g_2(x) = 3 \left( \frac{1-15x}{4} \right) \ln\left( \frac{1-15x}{4} \right)   $$
and $ \tau _k=\left\{ 0,\,\,0.009,\,\,0.014,\,\,0.022,\,\,0.03,\,\,0 04\right\}  $
meet the requirement. The choice of $ \tau _k $ is obviously not unique;
$ \displaystyle \left\{ 0.005k  \right\} _{k=1}^8 $
also works.

In checking the inequalities $ g_2(\tau _{i+1})\leq g_1(\tau _{i}) $, I merely use the 
values of the two sides, computed using some numerical computer software. For example,
MAPLE 13 gives
$$  g_2(0.04)= -0.8059612084... < -0.7882434741... = g_1(0.023).  $$
Some may object to this step on the ground that computer numerics contain 
rounding errors, thus making it unreliable. In fact, modern numerical
software is highly accurate with calculations involving well-known 
functions. As long as we keep ourselves within a reasonable scope of the accuracy
of the software, we should be very safe. MAPLE, by default,
does its calculations with an accuracy of 10 digits. As long as the
numerical inequality holds within a 5 or 6 digit accuracy, it should be
very reliable. In case of doubt, one can instruct MAPLE to use a
higher number of digits to recheck the inequality. For the extremely 
skeptical, they can employ interval arithmetic (available, for instance,
 as a toolbox of MATLAB) to quench any remaining doubt.

\newpage
The DIF technique is summarized as follows:

\begin{center}
\fbox{
\begin{minipage}{5in}
\par\vspace*{1mm}\par
Steps to prove a desired inequality \hspace*{20mm} {\Cr $f(x) > 0, \quad  x\in[\alpha ,\beta ]$}.
\begin{enumerate}
\item Find two $ {\Cr\ic} $ (${\Cr\dc}$) functions such that \qquad  {\Cr $ f(x)=g_1(x)-g_2(x) $.}
\item Find points: \hspace*{8mm}
$ {\Cr  \tau _1 < \tau _2 < \cdots < \tau _n } \hspace*{8.5mm} $
such that (\ref{g12}) holds.
\par\vspace*{1.0mm}\par
\end{enumerate}
\end{minipage}
}
\end{center}

To establish inequality (\ref{b4p}), one can use
$$  g_1(s) = \frac{\ln^3(3)}{3^s}  +  \frac{\ln^3(4)}{4^s} , \quad  g_2(s) =    \frac{\ln^3(2)}{2^s} + \frac{\ln^3(5)}{2\cdot5^s} ,   $$
which are obviously decreasing functions, and the uniformly spaced sequence
$$  \tau _k :  0,\,\, 0.02,\,\, ... \mbox{ (increment by 0.02)} ,\,\, 1 .  $$

Alternatively, one can use
$$  g_1(s) = - \frac{\ln^3(2)}{2^s} + \frac{\ln^3(3)}{3^s} \,, \qquad  g_2(s) =    \frac{\ln^3(4)}{4^s} - \frac{\ln^3(5)}{2\cdot5^s}  \, ,  $$
for which a much shorter sequence is sufficient:
$$  \tau _k : \quad  0,\,\, 0.4,\,\, 0.65,\,\, 0.8,\,\, 0.9,\,\, 1 .  $$
To complete the proof, one must fill in a necessary step, namely, to show that 
$ g_1(s) $ and $ g_2(s) $ are decreasing. To that end, we can try to show that
$$  g_1'(s) =  \frac{\ln^4(2)}{2^s} - \frac{\ln^4(3)}{3^s} < 0  $$
(and the analog $ g_2'(s)<0 $), for which the DIF technique is applicable. 

\par\vspace*{3mm}\par
\begin{REM} \em
The procedure has its limitations. It cannot handle non-strict inequalities, when
$ f(x)=0 $ at $ x=\alpha  $ or $ \beta  $ or some difficult-to-compute interior point. 
Sometimes it can serve as part of the proof, using another trick 
to take care of the non-strict inequality in a neighborhood of the critical point.
A similar remark applies to inequalities in a non-compact  interval such as $ [0,\infty ) $.
\end{REM}

\begin{REM} \em
If the inequality involves other variables, e.g.\ in 
the form of parameters, the procedure fails. Again, it may be possible to
find a way around the difficulty with some tricks.
\end{REM}

\begin{REM} \em
What can we do if $ g_1(x) $ and{/}or $ g_2(x) $ are not monotone?

If their critical points are known, we can simply verify the inequality in each
of the sub-intervals sub-divided by the critical points. To handle the general
case, we can choose a sufficient
large (positive or negative) $ n $ so that $ (x-\alpha +1)^n(g_1(x)-g_2(x)) $ is a difference
of two monotone functions. In particular, we can use this technique to verify if
a given value is a strict lower{/}upper bound of a given function.
\end{REM}

\begin{REM} \em
In practice, Step 1 together with a decent graphing of the functions $ g_1(s) $
and $ g_2(s) $, especially if the software allows you to zoom into any 
sub-interval of $ [\alpha ,\beta ] $, is sufficient to offer 99\% confidence that the desired
inequality is valid. The provision of the sequence $ \tau _k $ is mainly for the purpose
of allowing other people to double check the assertion, especially
in cases when $ g_1(s) $ and $ g_2(s) $ are extremely close together.
\end{REM}

\begin{REM} \em
For special classes of functions, 
there are other more efficient numerical procedures. For example,
the Sturm procedure is the best choice for
studying algebraic and trigonometric polynomials. See \cite{ak2}--\cite{k2} 
for details.
\end{REM}

\setcounter{equation}{0}\section{MAPLE Procedure for Determining{/}Verifying $ \tau _k $}

Finding or verifying $ \tau _k $, though conceptually straightforward, is
tedious if done manually. We present a simple MAPLE procedure that can
help in this aspect. 

Let us first describe the usage of the procedure. The procedure is
written in the MAPLE programming language and the code lines are typed
into a file, say called {\tt dif}. I used MAPLE 13 for the experiments.
Very likely, the same code works all for other versions.

\subsection*{Usage:}

\begin{center}
\fbox{
\small 
\begin{minipage}{6.5 in}
\begin{tabbing}\hspace*{0ex} \= \hspace*{30ex} \=\\
 \>  \hspace*{3mm} FIRST FORMAT \\
 \> {\tt dif(g1,g2,[a,b],...)}  \> {\tt g1}, {\tt g2} are expressions representing the functions $ g_1 $\\
 \>   \>   and $ g_2 $. $ a $ and $ b $ are the endpoints of the interval. \\
 \>   \>  In this form, the procedure outputs a sequence of \\
 \>   \>  numbers $ \tau _k $ satisfying (\ref{g12}). \\
 \>   \> \\
 \>   \> Additional options: \\
 \>   \>  \,\, {\tt long=1} \hspace*{5mm} Verbose output \\
 \>   \>  \,\, {\tt steps=N} \hspace*{3.3mm} Maximal number of $ \tau _k $ (default: 100) \\
 \>   \>  \,\, {\tt digits=n} \hspace*{1.5mm}  How many decimal points \\
 \>   \> \\
 \>  \hspace*{3mm} SECOND FORMAT \\
 \> {\tt dif(g1,g2,[t1,...,tn])}  \> Verifies if {\tt t1, ..., t1} satisfy (\ref{g12}) \\
\end{tabbing}
\end{minipage}
}
\end{center}

In an interactive MAPLE session, commands can be entered
either in Text or Math mode. Below, the command are shown in the Text mode
because it is easier to typeset. The Math mode
input and output are also shown.

While inside a MAPLE session, the
procedure is loaded by using the command 
{\tt read(dif)}.

Newer versions of MAPLE do not require the semicolon at the end of the
command if typed in the Math mode. Some older versions may need it.

\subsection*{Examples:}

\par\vspace*{\baselineskip}\par
{\small 
\begin{quote}
\begin{Verbatim}[frame=single]
read(dif);
g1 := ((1-3*x)/2)*ln((1-3*x)/2) + 2*((1-24*x)/5)*ln((1-24*x)/5);
g2 := 3*((1-15*x)/4)*ln((1-15*x)/4);
dif(g2, g1, [0, 0.04]);                   # FIRST FORMAT
dif(g1, g2, [0, 0.04]);
dif(g1, g2, [0, 0.04], long=1);
\end{Verbatim}
\end{quote}
}

\par\vspace*{\baselineskip}\par
\noindent In this first format, the third input argument 
$ [0,0.04] $ is a list of two numbers, representing $ [\alpha ,\beta ] $.
The output is shown below:

\par\vspace*{8mm}\par
\FG{140}{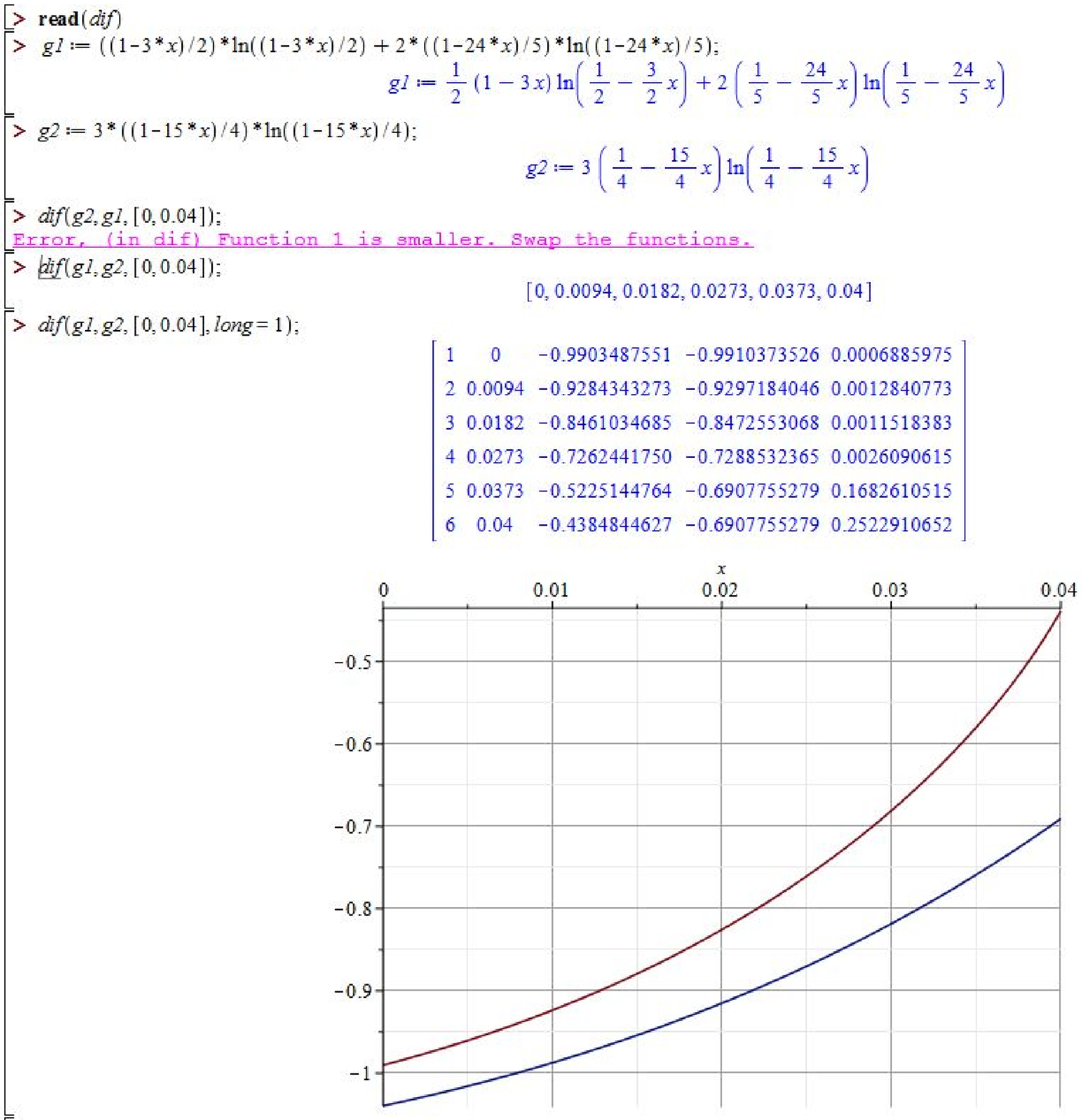}

\newpage
The procedure is written to mandate that the first argument
is the larger function. In case of an error, the
procedure will complain as shown. The procedure could have been easily
rewritten to automatically swap the functions in that situation. However,
I believe that the error checking feature serves as a further insurance
of the proper choice of {\tt g1} and~{\tt g2}.

In the {\tt long=1} version, the output is a matrix of five columns. The
first column is a count, namely $ k $, of the points $ \tau _k $ given in the
second column. If $ g_i(x) $ are increasing, the third column is $ g_1(\tau _k) $ and
the fourth column is $ g_2(\tau _{k+1}) $. If $ g_i(x) $ are decreasing, they are
$ g_1(\tau _{k+1}) $ and $ g_2(\tau _k) $ instead. In all cases, the values in the
third column should be larger than that in the same row of the fourth
column; the fifth column is the difference between these values.

The last three columns are listed as a visual check to make sure that 
the claimed inequalities (\ref{g12}) are valid within reasonable computational errors.
For example, if the values are computed with 10 (the default in MAPLE)
digit accuracy, then the differences between the values in the third and
fourth columns (i.e. the value in the fifth column) should be at least
$ 0.0001 $. If this is not true, the procedure can be rerun with increased
accuracy, for example, with the MAPLE command ``{\tt Digits := 20}''.

It is the user's responsibility to make sure that {\tt g1} and {\tt g2}
are either both increasing or both decreasing in $ [a,b] $. The procedure does not have
the intelligence to verify that premise rigorously.
In case that assumption is not valid, the procedure may still output 
some (useless) answer.

To help a bit, the {\tt long=1} option will plot the graphs of {\tt g1}
(red curve) and {\tt g2} (blue) to allow the user do a visual check.

If the monotonicity of $ g_i(x) $ is not obvious, try the DIF technique on $ g_i'(x) $.

\subsection*{Comments}
\begin{itemize}
\item The variable used in {\tt g1} and {\tt g2} can be arbitrary. It does not
matter if you have used $ x $ or $ y $; or $ s $.
The procedure automatically detects it.

\item If there is more than
one variable present in {\tt g1} and {\tt g2} (perhaps due to typos),
an error message will be printed. 
More complicated error checking is not available. It is the responsibility
of the user.

\item {\tt dif} fails if it takes more than 100 steps without reaching the 
right endpoint $ \beta  $ of the interval. This can mean either the desired inequality is
false or the difference between $ g_1(x) $ and $ g_2(x) $ so small that it 
actually requires more than 100 points of $ \tau _k $. One can try to exclude
the latter by using the {\tt steps=N} option with a larger {\tt N}.
\end{itemize}

\newpage
\subsection*{Further Examples}

By default, {\tt dif} rounds each $ \tau _k $ down to a number of decimal digits
compatible with about one-hundredth of the length of the interval $ [\alpha ,\beta ] $.
This can be altered with the option {\tt digits=n}. Using a higher {\tt n}
may generate a shorter $ \tau _k $ sequence. The specification of {\tt n} is only
suggestive. {\tt dif} may ignore it if necessary due to computational
accuracy.

The next two examples show the second format of the procedure. The 
third input argument  is a list of more than two numbers, which are interpreted
as a sequence of $ \tau _k $ to be verified.

\par\vspace*{\baselineskip}\par

{\small 
\begin{quote}
\begin{Verbatim}[frame=single]
dif(g1, g2, [0, 0.04], digits=5);
dif(g1, g2, [0, 0.009, 0.019, 0.025, 0.034, 0.04])        # SECOND FORMAT
dif(g1, g2, [0, 0.009, 0.017, 0.025, 0.034, 0.04])
\end{Verbatim}
\end{quote}
}

\par\vspace*{5mm}\par
The output is shown below.
\par\vspace*{5mm}\par

\FG{140}{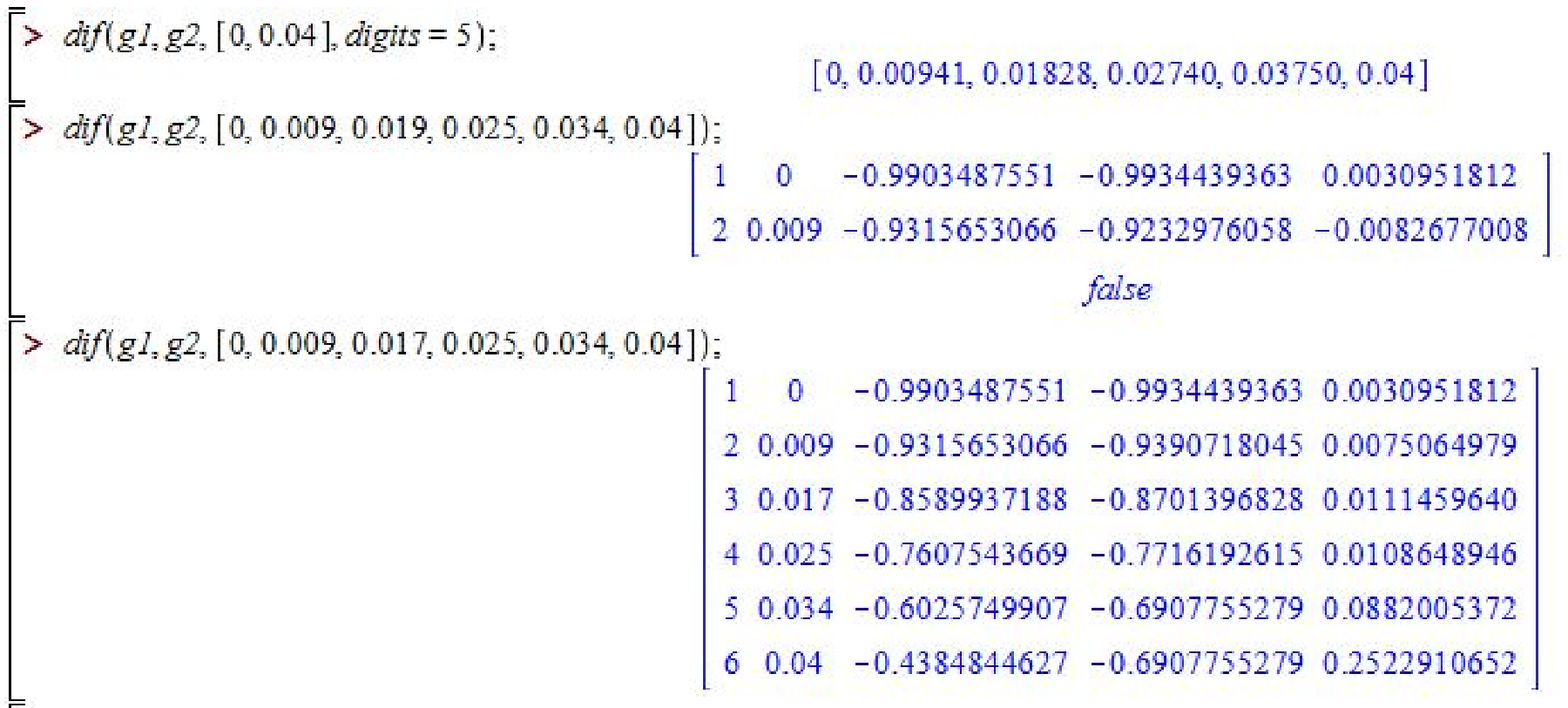}

\par\vspace*{\baselineskip}\par
When the supplied sequence $ \tau _k $ does not satisfy (\ref{g12}), the procedure will balk 
at the first instance of violation and print ``{\tt false}''. Failure
is also indicated by the fact that the last element of the last column
of the displayed matrix is negative.

In case of success,
the display is similar to the regular {\tt long=1} option, but without
the graphs of $ g_i(x) $. If graphs are desired, add {\tt long=1}.

\newpage
\setcounter{equation}{0}\section{Explanation of the MAPLE Code}

We briefly explain the code, assuming that the reader has rudimentary
knowledge of MAPLE programming.

The entire code (with added comments) is repeated in the Appendix. 
The reader can use cut and
paste to transfer that to a file, called ``{\tt dif}'',
on the computer.

We break the code (minus the comments)
up into small chunks.
Line numbers are added for reference and are not part of the code;
they should not be typed in the file. Blank lines have been added 
for aesthetic purposes and are optional.

\par\vspace*{4mm}\par
{ \footnotesize 
\begin{minipage}{5.5in}
\begin{Verbatim}[frame=single]
  1  ffloor := (x,n) \rightarrow  subsop(2=-n, convert(floor(x*10^n), float)):
\end{Verbatim}
\end{minipage}
}

\par\vspace*{1mm}\par
The first non-comment line in the file defines a separate 
procedure {\tt ffloor} which
is used by {\tt dif}, in line 27. It is a useful one to add to the
regular MAPLE repertoire. MAPLE has a built-in command {\tt floor(x)} to
round a floating point number {\tt x} down to the nearest integer, but there
is no command to conveniently round {\tt x} to the nearest float with
a given number of decimal places. One can achieve the effect by using
the command: \verb+floor(10^n*x)/10^n+.

The resulting value is correct. However, it is aesthetically undesirable because,
in most cases, it will be printed with some annoying trailing 0's.
Using {\tt ffloor} fixes that. 

The rest of the code is {\tt dif} proper.

\par\vspace*{4mm}\par
{ \footnotesize 
\begin{minipage}{5.5in}
\begin{Verbatim}[frame=single]
  2  dif := proc(g1,g2,A,\{long:=0,steps:=100,digits:=-100,relax:=99\})
  3    local aa, bb, ga, gb, f1, f2, tt, t1, cc, acc, tau, inc, var, ev;
\end{Verbatim}
\end{minipage}
}

\par\vspace*{1mm}\par
These two lines specify the procedure name {\tt dif}, the mandatory
arguments {\tt g1}, {\tt g2}, and {\tt A}, the optional arguments, and
the local variables.

\par\vspace*{4mm}\par
{ \footnotesize 
\begin{minipage}{5.5in}
\begin{Verbatim}[frame=single]
  4    ev := (a,g) \rightarrow  evalf(subs(var=a,g));
\end{Verbatim}
\end{minipage}
}

\par\vspace*{1mm}\par
The next line defines a procedure {\tt ev(a,g)} to be used only internally
within {\tt dif}, essentially as a shorthand for a frequent construction. 
It takes an expression {\tt g} and computes its numerical value when a 
float {\tt a} is substituted into the variable contained in {\tt g}.

\par\vspace*{4mm}\par
{ \footnotesize 
\begin{minipage}{5.5in}
\begin{Verbatim}[frame=single]
  5    var := indets(evalf([g1,g2]), name);
  6    if nops(var) > 1 then
  7       error "Functions have more than one variable: ", var;
  8    else
  9       var := var[1];
 10    end if;
\end{Verbatim}
\end{minipage}
}

\par\vspace*{1mm}\par
In line 5, {\tt indets(...)} returns a list of variables found in {\tt g1} and 
{\tt g2}, and stores it in {\tt var}.
Here {\tt evalf} is called before invoking
{\tt indets}, because the latter thinks $ \pi  $,
such as in $ \sin(2\pi x) $, as a variable 
instead of as a constant. If {\tt var} has more than one
element, then an error is thrown; otherwise, we change {\tt var} to
mean the sole variable it contains.

\newpage
\par\vspace*{-9mm}\par
{ \footnotesize 
\begin{minipage}{5.5in}
\begin{Verbatim}[frame=single]
 11    aa := A[1]; bb := A[-1]; 
 12    ga := ev(aa,g1); gb := ev(aa,g2);
 13    if ga < gb then error "Function 1 is smaller. Swap the functions."
 14    end if;
\end{Verbatim}
\end{minipage}
}

\par\vspace*{1mm}\par
Here we compute {\tt ga\! =} $ g_1(\alpha ) $ and {\tt gb\! =} $ g_2(\alpha ) $ and make sure that the former is 
larger.

\par\vspace*{4mm}\par
{\footnotesize 
\begin{minipage}{5.5in}
\begin{Verbatim}[frame=single]
 15    gb := ev(bb,g2);
 
 16    if ga < gb then f1 := g1; f2 := g2; inc := 1; 
 17    else f1 := -g2; f2 := -g1; inc := 0;
 18    end if;
\end{Verbatim}
\end{minipage}
}

\par\vspace*{1mm}\par
Next we compute {\tt gb\! =} $ g_1(\beta ) $, using it to decide whether $ g_i(x) $
are increasing or decreasing; this information is stored in the flag {\tt inc}.
Then {\tt f1} and {\tt f2} are assigned the appropriate
expression, to be used later to compute $ \tau _k $.

\par\vspace*{4mm}\par
{ \footnotesize 
\begin{minipage}{5.5in}
\begin{Verbatim}[frame=single]
 19    if nops(A) = 2 then 
          ...
 44    else
          ...
 58    end if;
 59    end proc:
\end{Verbatim}
\end{minipage}
}

\par\vspace*{1mm}\par
The rest of the code is divided into two separate cases: lines 20--42
when the user asks for $ \tau _k $, and lines 44--56 when the user wants to
verify if the supplied $ \tau _k $ work. Line 59 is the last line of the 
procedure and also of the file.

\par\vspace*{4mm}\par
{ \footnotesize 
\begin{minipage}{5.5in}
\begin{Verbatim}[frame=single]
 20       if digits > 0 then 
 21          acc := digits 
 22       else acc := 2-floor(log10(bb-aa))
 23       end if;
\end{Verbatim}
\end{minipage}
}

\par\vspace*{1mm}\par
These lines assign to {\tt acc} the number of digits.

\par\vspace*{4mm}\par
{ \footnotesize 
\begin{minipage}{5.5in}
\begin{Verbatim}[frame=single]
 24       tau := [aa];
 25       tt := aa; ga := ev(aa,f1); gb := ev(bb,f2); cc := 0;

 26       for cc from 1 to steps while (ga < gb) do
 27         t1 := min(ffloor((relax*fsolve(f2=ga,var=tt..bb)+tt)/
                  (relax+1),acc),bb);
 28         if t1 \leq  tt then acc := acc+1; next; end if;
 29         tt := t1; tau := [op(tau),tt];
 30         ga := ev(tt,f1);
 31       end do;
\end{Verbatim}
\end{minipage}
}

\par\vspace*{1mm}\par
The above is the loop to compute $ \tau _k $. Ideally, if numerical 
computation were exact (and all inequalities were replaced by equality), then
$ \tau _k $ would have been determined recursively by: 
$ \displaystyle \tau _k = g_2^{-1}(g_1(\tau _{k-1})). $
Geometrically, $ \tau _k $ is the successive values of the iterative algorithm
depicted in Figure 1.

\begin{figure}[h]
\begin{center}
\FG{120}{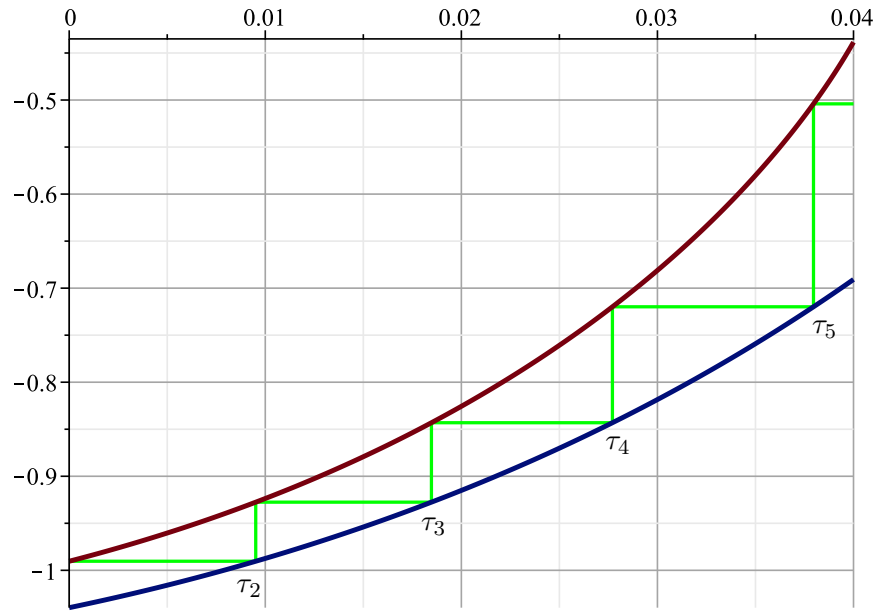} 
\par\vspace*{-1mm}\par
\caption{Iterative determination of $ \tau _k $.}
\end{center}
\par\vspace*{-20.5mm}\par \hspace*{48.5mm} $ \tau _2 $
\par\vspace*{-13mm}\par \hspace*{73mm} $ \tau _3 $
\par\vspace*{-15.5mm}\par \hspace*{97.5mm} $ \tau _4 $
\par\vspace*{-20.3mm}\par \hspace*{125.0mm} $ \tau _5 $
\end{figure}

\par\vspace*{40mm}\par

In practice, we have to make $ \tau _k $ smaller; this is done
in line 27 using {\tt relax}. By default {\tt relax=99}, and that can
be altered by the user. Choosing a smaller {\tt relax} makes $ \tau _k $
smaller. 

The variable {\tt cc} counts the number of 
iterations. The loop exits either if {\tt cc} exceeds the assigned number 
of steps or if $ \tau _k $ has reached the right endpoint
$ \beta  $ of the interval. 

\par\vspace*{4mm}\par
{ \footnotesize 
\begin{minipage}{5.5in}
\begin{Verbatim}[frame=single]
 32       if ga < gb then 
 33          cc := min(steps,5);
 34          printf("    **** Fails after %d steps. Last %d steps are: ",
                steps, cc); 
 35          print(tau[-cc..-1]);
\end{Verbatim}
\end{minipage}
}

\par\vspace*{1mm}\par
In the former case, a failure message is displayed.

\par\vspace*{4mm}\par
{ \footnotesize 
\begin{minipage}{5.5in}
\begin{Verbatim}[frame=single]
 36       else
 37          tau := [op(tau),bb]; 
 38          if long = 1 then 
 39             dif(g1,g2,tau); 
 40             plot([g1,g2],var=aa..bb,gridlines=true,size=[550,400]);
 41          else tau; 
 43          end if;
 44       end if;
\end{Verbatim}
\end{minipage}
}

\par\vspace*{1mm}\par
Otherwise, line 36 adds $ \beta  $ as the last $ \tau _k $. If the {\tt long=1} option
is set, {\tt dif} is invoked with the second format to print the verbose
output, and the graphs of $ g_i(x) $ are plotted. Without the verbose option,
only the sequence {\tt tau} is displayed.

\newpage
{ \footnotesize 
\begin{minipage}{5.5in}
\begin{Verbatim}[frame=single]
 45       t1 := nops(A);
 46       for cc from 1 to t1-1 do
 47          ga := ev(A[cc+1-inc],g1); gb := ev(A[cc+inc],g2);
 48          if cc = 1 then
 49             tau := <1 | A[1] | ga | gb | ga-gb >;
 50          else
 51             tau := <tau,<cc | A[cc] | ga | gb | ga-gb >>;
 52          end if;
 53          if ga < gb then print(tau); return false; end if;
 54       end do;
 55       ga := ev(A[-1],g1); gb := ev(A[-1],g2);
 56       tau := <tau,<t1 | A[-1] | ga | gb | ga-gb >>;
 57       print(tau);
\end{Verbatim}
\end{minipage}
}

\par\vspace*{1mm}\par
This final part of the code handles the {\tt long} display when {\tt dif}
is invoked in the second format, with supplied $ \tau _k $. In this
format, $ \tau _k $ is given in the input variable {\tt A}. The variable {\tt tau},
on the other hand, is
free to be used for another purpose, namely, to store the matrix 
containing the five columns of numbers to be displayed at the end.
Lines 45--53 is the loop to build the rows (except the last one
which is filled in by line 55) of the matrix. 
Line 53 prints the partial {\tt tau} and returns ``{\tt false}'' 
in case condition (\ref{g12}) fails.
Line 56 prints the matrix
and the job is done.

\newpage
\section*{Appendix: The Entire MAPLE Code}

\begin{quote}
{ \footnotesize 
\begin{verbatim}
# CUT AND PAST THE FOLLOWING LINES INTO A FILE CALLED "dif"
 
#  ffloor(x,n)         ROUND FLOAT x DOWN TO n DECIMAL PLACES
ffloor := (x,n) -> subsop(2=-n, convert(floor(x*10^n), float)):

#  dif(g1,g2,A,...)    SEE SECTION ON USAGE
dif := proc(g1,g2,A,{long:=0,steps:=100,digits:=-100,relax:=99})
  local aa, bb, ga, gb, f1, f2, tt, t1, cc, acc, tau, inc, var, ev;

  #  ev(a,g)           EVALUATE g TO FLOAT, SUBSTITUTING var=a
  ev := (a,g) -> evalf(subs(var=a,g));

  # DETERMINES THE VARIABLE IN g
  var := indets(evalf([g1,g2]), name);
  if nops(var) > 1 then
     error "Functions have more than one variable: ", var;
  else
     var := var[1];
  end if;

  # CHECK IF g1 > g2
  aa := A[1]; bb := A[-1]; 
  ga := ev(aa,g1); gb := ev(aa,g2);
  if ga < gb then error "Function 1 is smaller. Swap the functions."
  end if;

  # DETERMINE IF g1 IS INCREASING OR DECREASING
  gb := ev(bb,g1);
  if ga < gb then f1 := g1; f2 := g2; inc := 1; 
  else f1 := -g2; f2 := -g1; inc := 0;
  end if;

  if nops(A) = 2 then 

     # FIRST FORMAT OF INVOCATION: DETERMINE [ tau_k ]
     if digits > 0 then 
        acc := digits 
     else acc := 2-floor(log10(bb-aa))
     end if;
     tau := [aa];
     tt := aa; ga := ev(aa,f1); gb := ev(bb,f2); cc := 0;

     # LOOP TO DETERMINE tau_k
     for cc from 1 to steps while (ga < gb) do
       t1 := min(ffloor((relax*fsolve(f2=ga,var=tt..bb)+tt)/(relax+1),acc),bb);
       if t1 <= tt then acc := acc+1; next; end if;
       tt := t1; tau := [op(tau),tt];
       ga := ev(tt,f1);
     end do;

     if ga < gb then 
        # FAIL
        cc := min(steps,5);
        printf("    **** Fails after %d steps. Last %d steps are: ", steps, cc);
        print(tau[-cc..-1]);
     else
        # SUCCESS
        tau := [op(tau),bb]; 
        if long = 1 then 
           dif(g1,g2,tau); 
           plot([g1,g2],var=aa..bb,gridlines=true,size=[550,400]);
        else tau; 
        end if;
     end if;

  else

     # SECOND FORMAT OF INVOCATION: VERIFY [ tau_k ]
     t1 := nops(A);

     # BUILD MATRIX tau
     for cc from 1 to t1-1 do
        ga := ev(A[cc+1-inc],g1); gb := ev(A[cc+inc],g2);
        if cc = 1 then
           tau := <1 | A[1] | ga | gb | ga-gb >;
        else
           tau := <tau,<cc | A[cc] | ga | gb | ga-gb >>;
        end if;
        if ga < gb then print(tau); return false; end if;
     end do;

     ga := ev(A[-1],g1); gb := ev(A[-1],g2);
     tau := <tau,<t1 | A[-1] | ga | gb | ga-gb >>;

     print(tau);
  end if;
  end proc:                         # END OF FILE
\end{verbatim}
}
\end{quote}
\end{document}